\documentclass[11pt]{article}

\usepackage{amsthm} 
\usepackage{fullwidth}
\usepackage[top=1.5in, bottom=1.5in, left=1in, right=1in]{geometry}
\usepackage{xypic}
\usepackage{amsgen}
\usepackage{amsmath}
\usepackage{amstext}
\usepackage{amsbsy}
\usepackage{amsopn}
\usepackage{amsfonts}
\usepackage{amssymb}
\usepackage{eepic}
\usepackage[pdftex]{graphicx}
\usepackage{epsf}
\usepackage[usenames,dvipsnames]{xcolor}
\usepackage{pstricks}
\usepackage{makeidx}
\xyoption{all}

\usepackage{soul} 
\definecolor{aqua}{rgb}{0.0, 1.0, 1.0}

\definecolor{eggshell}{rgb}{0.94, 0.92, 0.84}
\definecolor{cosmiclatte}{rgb}{1.0, 0.97, 0.91}
\definecolor{cornsilk}{rgb}{1.0, 0.97, 0.86}
\definecolor{blond}{rgb}{0.98, 0.94, 0.75}
\definecolor{beige}{rgb}{0.96, 0.96, 0.86}
\definecolor{lemonchiffon}{rgb}{1.0, 0.98, 0.8}
\definecolor{pastelyellow}{rgb}{0.99, 0.99, 0.59}
\definecolor{unmellowyellow}{rgb}{1.0, 1.0, 0.4}

\usepackage{enumerate}

\def\Box{\square}

\def\tra#1{\smash{\mathop{\mid\kern
-1pt\joinrel\relbar\joinrel\relbar}\limits^{*}_{#1}}}
\def\longtra#1{\smash{\mathop{\mid\kern
-1pt\joinrel\relbar\joinrel\relbar\joinrel\relbar}\limits^{*}_{#1}}}
\def\vlongtra#1{\smash{\mathop{\mid\kern-1pt\joinrel\relbar\joinrel\relbar\joinrel\relbar\joinrel\relbar}\limits^{*}_{#1}}}
\def\vvlongtra#1{\smash{\mathop{\mid\kern
-1pt\joinrel\relbar\joinrel\relbar\joinrel\relbar\joinrel\relbar\joinrel\relbar}\limits^{*}_{#1}}}
\def\vvvlongtra#1{\smash{\mathop{\mid\kern
-1pt\joinrel\relbar\joinrel\relbar\joinrel\relbar\joinrel\relbar\joinrel\relbar\joinrel\relbar}\limits^{*}_{#1}}}
\def\etra#1{\smash{\mathop{\mid\kern
-1pt\joinrel\relbar\joinrel\relbar}\limits_{#1}}}

\def\oo{\overline}

\def\rk{{\rm rk}}

\def\min{\mbox{min}\,}

\def\N{\mathbb{N}}

\def\ZZ{\mathbb{Z}}

\def\p{\varphi}

\def\inv{^{-1}}
\def\la{\langle}
\def\ra{\rangle}

\def\bV{{\bf V}} 
 
\def\bAb{{\bf Ab}} 
\def\tF{{\text{F}}}
\def\tFA{{\text{FA}}}
\def\tFS{{\text{FS}}}

\newcommand{\rkG}{{\text{rk}_\text{G}}}
\newcommand{\rkM}{{\text{rk}_\text{M}}}
\newcommand{\rkS}{{\text{rk}_\text{S}}}

\renewcommand{\geq}{\geqslant}
\renewcommand{\leq}{\leqslant}

\def\bi{\begin{itemize}}
\def\ei{\end{itemize}}
\def\ben{\begin{enumerate}}
\def\een{\end{enumerate}}
\def\beq{\begin{equation}}
\def\eeq{\end{equation}}

\newtheorem{T}{Theorem}[section]
\newcommand{\bt}{\begin{T}}
\newcommand{\et}{\end{T}}
\newcommand{\ftd}{$\square$\end{T}}

\newtheorem{Proposition}[T]{Proposition}
\newcommand{\bp}{\begin{Proposition}}
\newcommand{\ep}{\end{Proposition}}
\newcommand{\fpd}{$\square$\end{Proposition}}

\newtheorem{Definition}[T]{Definition}
\newcommand{\bd}{\begin{Definition}}
\newcommand{\ed}{\end{Definition}}

\newtheorem{Lemma}[T]{Lemma}
\newcommand{\bl}{\begin{Lemma}}
\newcommand{\el}{\end{Lemma}}
\newcommand{\fld}{$\square$\end{Lemma}}

\newtheorem{Corol}[T]{Corollary}
\newcommand{\bc}{\begin{Corol}}
\newcommand{\ec}{\end{Corol}}
\newcommand{\fcd}{$\square$\end{Corol}}

\newtheorem{Result}[T]{Result}
\newcommand{\br}{\begin{Result}}
\newcommand{\er}{\end{Result}}
\newcommand{\frd}{$\square$\end{Result}}

\newtheorem{Remark}[T]{Remark}
\newcommand{\brem}{\begin{Remark}}
\newcommand{\erem}{\end{Remark}}
\newcommand{\fremd}{$\square$\end{Remark}}

\newtheorem{Example}[T]{Example}
\newcommand{\be}{\begin{Example}}
\newcommand{\ee}{\end{Example}}

\newtheorem{Problem}[T]{Problem}
\newcommand{\bq}{\begin{Problem}}
\newcommand{\eq}{\end{Problem}}

\def\qed{$\Box$
\par\bigbreak}

\def\proof
  {\par\medbreak\noindent{\bf Proof}.\enspace}

\newlength{\lengtha} 
\newlength{\lengthb} 

\def\abstract#1{\par\bigskip
\begingroup\small
\baselineskip=12truept
\begin{center}ABSTRACT\end{center}
\par\medskip\par\noindent
\null\hfill\hbox{\vbox{\hsize=5truein\noindent#1}}
\hfill\null\par\endgroup\par}

\title{On the semigroup rank of a group}
\author{{\bf M\'ario J.\ J.\ Branco, Gracinda M.\ S.\ Gomes and Pedro V.\ Silva}}

\date{\today}

\begin{document}

\maketitle

\begin{center}\small

\bigskip

\end{center}

\abstract{%
For an arbitrary group $G$, it is shown that
either the semigroup rank $G\rkS$ equals the group rank $G\rkG$, or $G\rkS = G\rkG+1$. 
This is the starting point for the rest of the article, 
where the semigroup rank for diverse kinds of groups is analysed. 
The semigroup rank of relatively free groups, for any variety of groups, is computed. 
For a finitely generated abelian group~$G$, it is proven that 
$G\rkS = G\rkG+1$ if and only if $G$ is torsion-free.
In general, this is not true. Partial results are obtained in the nilpotent case. 
It is also proven that if $M$ is a connected closed surface, then 
$(\pi_1(M))\rkS = (\pi_1(M))\rkG+1$ if and only if $M$ is orientable.
%
}


\section{Introduction} \label{Sec:intro}


Given a semigroup $S$ and $X \subseteq S$, we denote by $X^+$ the subsemigroup of $S$ 
generated  by $X$. 
If $S$ is a monoid (respectively a group), 
we denote by $X^*$ (respectively $\langle X\rangle$) the submonoid 
(respectively subgroup) of $S$ generated by $X$. 
For a group~$G$ and a nonempty subset $X$ of $G$, one has 
$\langle X\ra = ( X \cup X\inv )^+$, where 
$X\inv = \{x\inv \colon x \in X \}$, 
and hence 
$G$ is finitely generated as a group 
if and only if $G$ is finitely generated as a semigroup. 
Thus, the expression ``finitely generated'' 
does not need any modifier to be unambiguous. 

The {\em semigroup rank} of a semigroup $S$ is defined as
$$S\rkS = \left\{
\begin{array}{ll}
\min\{ |X| \colon X \subseteq S\mbox{ and } X^+ = S\}&\mbox{ if $S$ is finitely generated}\\
\infty&\mbox{ otherwise}
\end{array}
\right.$$

Analogously, the {\em monoid rank} of a monoid $M$ is defined as
$$M\rkM = \left\{
\begin{array}{ll}
\min\{ |X| \colon X \subseteq M\mbox{ and } X^* = M\}&\mbox{ if $M$ is finitely generated}\\
\infty&\mbox{ otherwise}
\end{array}
\right.$$
and the {\em group rank} of a group $G$ is defined as
$$G\rkG = \left\{
\begin{array}{ll}
\min\{ |X| \colon X \subseteq G\mbox{ and } \la X\ra = G\}&\mbox{ if $G$ is finitely generated}\\
\infty&\mbox{ otherwise}
\end{array}
\right.$$

We shall refer to a semigroup (respectively monoid, group) generating set of minimum size as a 
{\em semigroup basis} (respectively {\em monoid basis}, {\em group basis}).

The rank of a semigroup has been widely studied, see for example~\cite{GH87, FQ11, Gray14, DEM16}.  
Most investigation on semigroup rank aims at finding  
the semigroup rank of specific finite semigroups, such as 
semigroups of transformations, Rees matrix semigroups, and endomorphism monoids of algebraic structures. 
In~\cite{Gray14} it is given a general theory that encompasses many results on semigroup rank. 

For a monoid and for a group, as mentioned we have different types of ranks,  
and so it is a natural question to compare $M\rkM$ and $M\rkS$ for a monoid~$M$, as well as $G\rkG$, $G\rkM$ and $G\rkS$ for a group~$G$. 
As one may guess, in the finite case this comparison is not of great interested, 
however the situation becomes rather complicated when one leaves the finite environment. 
Studying semigroup bases, and so the semigroup rank, of a group is certainly an important topic, 
in particular because semigroup bases for finitely generated infinite groups 
arise naturally from semigroup presentations involved in several subjects of Combinatorial Group Theory, 
usually related to algorithmic issues. 
Well-known instances appear in the theory of automatic groups~\cite{Eps} 
and in the use of rewriting systems as a tool~\cite{HEO}.

This article is essentially devoted to the comparison of $G\rkS$ and $G\rkG$ for a group $G$.   
We start by comparing all possible ranks for monoids and for groups in the second section. 
The difficulty is in writing the semigroup rank (or monoid rank) of a group in terms 
of its group rank.
Despite the fact that, in general, 
the semigroup rank is equal to the group rank or equal to the group rank plus one, 
to establish for a given group their equality may prove to be very hard.  
The subsequent sections concern groups, for which we analyse both 
the group and the semigroup ranks.  
We look first at relatively free groups, for which we completely determine the  
semigroup rank. 
Thereafter, we consider arbitrary groups with some specificity, 
namely with/without torsion, abelian, and more generally nilpotent. 
We prove that if $G$ is a finitely generated abelian group, 
then $G\rkS = G\rkG+1$ if and only if $G$ is torsion-free. 
We end with the analysis of surface groups.  
It is shown that if $M$~is a connected closed surface, 
then $(\pi_1(M))\rkS = (\pi_1(M))\rkG+1$ if and only if $M$ is orientable.

\section{First results}

To compare $M\rkM$ and $M\rkS$ for a monoid $M$ is quite easy. 
Indeed, if $M = X^+$, then $M \setminus \{ 1 \} \subseteq (X \setminus \{ 1 \})^+$,  
and so $1$ belongs to a semigroup basis of $M$ if and only if $1 \notin (M \setminus \{ 1 \})^+$, 
when 1 must indeed belong to any semigroup generating set of $M$. 
Therefore we obtain:

\bp
\label{sm}
Let $M$ be a monoid. Then:
$$M\rkS = \left\{
\begin{array}{ll}
M\rkM+1&\mbox{ if $M$ is finitely generated and }1 \notin (M \setminus \{ 1 \})^+\\
M\rkM&\mbox{ otherwise}
\end{array}
\right.$$
\ep

We immediately get:

\bc
\label{smg}
Let $G$ be a group.Then:
$$G\rkS = \left\{
\begin{array}{ll}
G\rkM+1&\mbox{ if $G$ is trivial}\\
G\rkM&\mbox{ otherwise}
\end{array}
\right.$$
\ec

By these results, from now on we will only focus on the comparison of $G\rkS$ and $G\rkG$ for a group $G$. 
Since, as we have noticed before,  
$G$ is finitely generated as a group if and only it is finitely generated as a semigroup, 
$G\rkS = G\rkG$ if $G$ is not finitely generated. 
We shall restrict henceforth our attention to finitely generated groups.

As every semigroup generating set of a group $G$ must also be a group generating set, we get the inequality 
$$G\rkG \leq G\rkS.$$
This inequality will be used throughout the text without further reference.

\bp \label{prop:GgGen-SmgGen}
If\/ $\{ a_1, \dotsc, a_n \}$ is a group generating set for 
a group $G$, then \linebreak 
$\bigl\{ a_1, \dotsc, a_n, a_n\inv\dotsm a_1\inv \bigr\}$ 
is a semigroup generating set for $G$.
\ep

\proof 
Let $G$ be a group, and let 
$A = \{ a_1, \dotsc, a_n \} \subseteq G$. 
Consider the element $a = a_n\inv\dotsm a_1\inv$. 
Since 
$$a_i\inv = a_{i+1}\dotsm a_{n}\bigl(a_n\inv\dotsm a_1\inv\bigr)a_1\dotsm a_{i-1}  
    = a_{i+1}\dotsm a_{n}\cdot a \cdot a_1\dotsm a_{i-1}$$
for $i = 1,\dotsc,n$, 
we have 
$\la A \ra = (A \cup A\inv)^+ = \{ a_1, \dotsc, a_n, a \}^+$, 
and hence the result follows. 
\qed

Proposition~\ref{prop:GgGen-SmgGen} has the following immediate consequence for the general case. 

\bc
\label{general}
Let $G$ be a finitely generated group. Then 
$G\rkS = G\rkG$ or $G\rkS = G\rkG+1$. 
\ec

\proof
We have seen that $G\rk_G \leq G\rk_S$. 
By Proposition~\ref{prop:GgGen-SmgGen}, $G\rk_S \leq G\rk_G +1$, 
whence $G\rkS = G\rkG$ or $G\rkS = G\rkG+1$.
\qed

This corollary does not say, however, which of the equalities, 
$G\rk_S = G\rk_G$ or $G\rk_S = G\rk_G +1$, 
holds for a given group~$G$. 
We answer this question for various kinds of groups in the next sections.

\section{Relatively free groups}

Given a variety $\bV$ of groups and a set $A$, we denote by $\tF_A(\bV)$ 
the (relatively) free group in~$\bV$ on~$A$. 
We shall assume, for simplicity, that V is nontrivial and $A \subseteq F_A(V)$.
It is well-known that, given sets $A$ and $B$, 
if $|A| = |B|$, then $\tF_A(\bV) \cong \tF_B(\bV)$. 
Thus, we follow the standard notation and denote 
$\tF_A(\bV)$, where $|A|=n$, by $\tF_n(\bV)$.
The (absolutely) free group (respectively free abelian group) on~$A$ 
will be denoted simply by~$\tF_A$ (respectively~$\tFA_A$), 
and the free group (respectively free abelian group) on a set 
of cardinality~$n$ will be denoted by
$\tF_n$ (respectively~$\tFA_n$). 
In the sequel we will use the fact that 
$(\tF_n(\bV))\rkG = n$, 
for every nontrivial variety $\bV$ of groups and $n \geq 0$  
(see~\mbox{\cite[13.53]{HNeu67}}).


A {\em group presentation} is a formal expression of the form 
Gp$\langle A \mid S\ra$, where $A$ is a set and $S$ is a subset of $\tF_A$, the free group on $A$. 
The group defined by this presentation is the quotient $\tF_A/\la\!\la S\ra\!\ra$, 
where $\la\!\la S\ra\!\ra$ denotes the normal subgroup of $\tF_A$ generated by the subset $S$.


As usual, we consider the (semi)groups defined by presentations up to isomorphism. 

A group $G$ is said to be {\em hopfian} if every surjective morphism $G \to G$ is necessarily an isomorphism.

We aim to determine the semigroup rank of $\tF_n(\bV)$ for all varieties $\bV$ of groups. 
For this purpose we start with the following.

\bl \label{prop:FreeV-Hopfian-RkIgualdade}
Let $\bV$ be a nontrivial variety of groups and let $A$ be a nonempty finite set 
such that $\tF_A(\bV)$~is hopfian. We have 
$(\tF_A(\bV))\rkS = (\tF_A(\bV))\rkG$ 
if and only if 
$\tF_A(\bV) = A^+$. 
\el 

\proof
Assume that $(\tF_A(\bV))\rkS = (\tF_A(\bV))\rkG$. 
Let us denote the free semigroup on the set~$A$ by~$\tFS_A$. 
Consider the morphism $\eta \colon \tFS_A \to \tF_A(\bV)$ 
such that $a\eta = a$ for any $a \in A$. 
Since $(\tF_A(\bV))\rkS = |A|$, there exists a surjective morphism 
$\mu \colon \tFS_A \to \tF_A(\bV)$. 
It follows that, by the properties of free groups, 
there exists a morphism 
$\varphi \colon \tF_A(\bV) \to \tF_A(\bV)$ 
such that $\eta_{|_A} \varphi = \mu_{|_A}$. 
Therefore $\eta \varphi = \mu$, and hence, since $\mu$ is surjective, $\varphi$ is also surjective. 
Then $\varphi$ is an isomorphism since $\tF_A(\bV)$ is hopfian.
Now, the fact that $\varphi$ is an isomorphism and $\mu$ is surjective 
implies that $\eta$~has to be surjective, and hence 
$\tF_A(\bV) = A^+$. 

Conversely, suppose that $\tF_A(\bV) = A^+$. 
Then $(\tF_A(\bV))\rkS \leq |A|$. 
But we also have $|A| = (\tF_A(\bV))\rkG \leq (\tF_A(\bV))\rkS$, 
whence 
$(\tF_A(\bV))\rkS = (\tF_A(\bV))\rkG$. 
\qed

\bp
\label{hi}
Let $G$ be a finitely generated group with $G\rkG = n$. 
If $\tFA_n$ is a homomorphic image of\/~$G$, then $G\rkS = G\rkG+1 = n+1$.
\ep

\proof
We start by showing that $(\tFA_n)\rkS = n+1$. 
Let $A$ be a set such that $|A| = n$, and consider 
$\tFA_n = \tF_A({\bAb})$, where $\bAb$ denotes the variety of abelian groups.  
Since $(\tFA_n)\rk_G = n$, by Corollary~\ref{general} either 
$(\tFA_n)\rk_S = n$ or $(\tFA_n)\rk_S = n+1$. 
Suppose that $(\tFA_n)\rkS = n$. 
By \cite[32.1 and 41.44]{HNeu67}, the group $\tFA_n$ is hopfian.
Then $\tFA_n = A^+$, by Lemma~\ref{prop:FreeV-Hopfian-RkIgualdade}. 
By the properties of the free abelian group, there exists a morphism 
$\varphi \colon \tFA_n \to \ZZ$ such that $a\varphi = 1$, 
for every $a \in A$. 
On the one hand, this morphism is surjective, and, 
on the other hand, we have 
$(\tFA_n)\varphi = A^+\varphi = (A\varphi)^+  = \{ 1 \}^+ = \N$, a contradiction. 
Therefore  $(\tFA_n)\rkS = n+1$. 

Now, assume that there exists a surjective morphism $\p \colon G \to \tFA_n$. 
Again by Corollary~\ref{general},  either 
$G\rk_S = G\rk_G = n$ or $G\rk_S = G\rk_G + 1 = n+1$. 
Suppose that $G\rkS = G\rkG$. 
Let $B \subseteq G$ be a semigroup basis of $G$. 
Then $B\p$ generates $\tFA_n$ as a semigroup, hence   
$(\tFA_n)\rkS \leq |B\p| \leq |B| = n$, a contradiction. 
Therefore $G\rkS = G\rkG+1$ as required.
\qed

In the next theorem, we compare $(\tF_A(\bV))\rkS$ with $(\tF_A(\bV))\rkG$ 
for an arbitrary variety $\bV$ of groups. 
Recall that, given a variety $\bV$ of groups, either 
$\bV$ satisfies some identity $x^k = 1$, where $k \geq 1$, or 
$\bV$~contains the additive group $\ZZ$. 
The former are called {\em periodic varieties}. 

\bt \label{prop:FreeV-Hopfian-rk}
Let $\bV$ be a nontrivial variety of groups and let $n \geq 1$. 
The following holds: \vspace{-1mm}
\begin{enumerate}[(i)]
 \item \label{prop:FreeV-Hopfian-rk-i} 
  If $\bV$ is periodic, then 
  $(\tF_n(\bV))\rkS = (\tF_n(\bV))\rkG = n$; \vspace{-1mm}
 \item \label{prop:FreeV-Hopfian-rk-ii} 
  If $\ZZ \in \bV$, then 
   $(\tF_n(\bV))\rkS = (\tF_n(\bV))\rkG + 1 = n+1$.  
\end{enumerate}
\et

\proof 
\eqref{prop:FreeV-Hopfian-rk-i} 
Assume that $\bV$ is periodic. 
Let $A$ be a group basis of $\tF_n(\bV)$. 
Then, for each $a \in A$, the element $a\inv$ is a power of $a$ of positive exponent, 
whence $a\inv \in A^+$. 
It follows that $\tF_n(\bV) = \langle A \rangle = (A \cup A\inv)^+ = A^+$, and hence 
$(\tF_n(\bV))\rkS \leq |A| = n$. 
Since we also have $n = (\tF_n(\bV))\rkG \leq (\tF_n(\bV))\rkS$, 
we obtain 
$(\tF_n(\bV))\rkS = (\tF_n(\bV))\rkG = n$. 

\eqref{prop:FreeV-Hopfian-rk-ii} 
Assume that $\ZZ \in \bV$. 
Thus the direct power $\ZZ^n$ belongs to $\bV$. 
Since $\tFA_n \cong \ZZ^n$, 
the group $\tFA_n$ is a homomorphic image of $\tF_n(\bV)$, 
and so, by Proposition~\ref{hi},  
$(\tF_n(\bV))\rkS = (\tF_n(\bV))\rkG + 1 = n+1$.~\qed

One of the consequences of Theorem~\ref{prop:FreeV-Hopfian-rk} is 
the following. 

\bc
\label{con}
For each $n \geq 1$ and any variety $\bV$ of groups, 
every group basis of~$\tF_n(\bV)$ is contained in some semigroup basis. 
In particular, every group basis of $\tF_n$ is contained in some semigroup basis.
\ec

\proof
Let $n \geq 1$ and let $\bV$ be a variety of groups. 
We saw in the proof of Theorem~\ref{prop:FreeV-Hopfian-rk} 
that, when $\bV$ is periodic, every group basis of $\tF_n(\bV)$ is also a semigroup basis. 

In case $\bV$ contains $\ZZ$, the conclusion is immediate by 
Theorem~\ref{prop:FreeV-Hopfian-rk} and Proposition~\ref{prop:GgGen-SmgGen}.  
\qed

Next we present counterexamples for two natural questions concerning groups: 
does every semigroup basis contain a group basis? 
Does every semigroup generating set contain a semigroup basis? 
The answer is negative for both. 
Bear in mind that $\mathbb{Z} \cong \tF_1$. 

\be
\label{ce1}
The set 
$\{2,-3\}$ is a semigroup basis of $\mathbb{Z}$, 
but contains no group basis.
\ee

\be
\label{ce2} 
The set 
$\{-6,2,3\}$ is a semigroup generating set of $\mathbb{Z}$ which contains no semigroup basis.
\ee

Finally, we consider decidability. 

\bp
\label{deci}
It is decidable whether or not a given subset of $\tF_n$ is a semigroup basis.
\ep

\proof
Let $S \subseteq \tF_n$. 
In view of Theorem~\ref{prop:FreeV-Hopfian-rk}\eqref{prop:FreeV-Hopfian-rk-ii} applied to $\tF_n$, 
we may assume that $|S| = n+1$. 
By Benois' Theorem~\cite{Benois87}, the language $\oo{S^+}$ constituted by the reduced forms of words in $S^+$ is an effectively constructible rational language, as it is $\oo{\tF_n}$. 
But $S$ is a semigroup basis of $\tF_n$ if and only if $\oo{S^+} = \oo{\tF_n}$, 
and this equality can be decided since equality is decidable for rational languages~\cite[Sec.\ 3.3]{HU79}.
\qed

\section{Arbitrary groups}

Since every semigroup basis of a group $G$ is necessarily a group generating set for~$G$, 
it follows that $G\rkS = G\rkG$ if and only if $G$ admits a group basis which is also a semigroup basis.
However it looks to be too hard to characterize precisely all the instances in which $G\rkS = G\rkG$. 
The discussion we initiate here, which depends in most cases on the property of being torsion-free or not, 
may shed some light on the range of possibilities to tackle this problem.

We start by introducing a sufficient condition to ensure $G\rkS = G\rkG$.

\bp
\label{torsion}
Let $G$ be a finitely generated group. 
If some group basis of\/ $G$ contains an element of finite order, then $G\rkS = G\rkG$.
\ep

\proof
Assume that $\{ a_1, \ldots,a_n\}$ is a group basis of $G$ and $a_n$ has order $m \geq 2$. 
By Proposition~\ref{prop:GgGen-SmgGen}, 
$\bigl\{ a_1, \dotsc, a_n, a_n\inv\dotsm a_1\inv \bigr\}$ 
is a semigroup generating set for $G$. 
We claim that \linebreak 
$S = \bigl\{ a_1,\ldots,a_{n-1}, a_n\inv\dotsm a_{1}\inv \bigr\}$ 
is a semigroup generating set for $G$. 
Indeed, 
$$a_n\inv = (a_n\inv\dotsm a_{1}\inv)a_1\dotsm a_{n-1} \in S^+$$
and so $a_n = (a_n\inv)^{m-1} \in S^+$ too. 
Then $G\rkS \leq G\rkG$. 
Since also $G\rkG \leq G\rkS$, we get equality as required.
\qed

We immediately get the following generalization of 
Theorem~\ref{prop:FreeV-Hopfian-rk}\eqref{prop:FreeV-Hopfian-rk-i}, 
as in that case $\tF_n(\bV)$ is a torsion group. 

\bc
\label{torsion2}
If\/ $G$ is a torsion group, then $G\rkS = G\rkG$.
\ec

Notice that the finitely generated torsion-free abelian groups are precisely the 
finitely generated free abelian groups~\cite[Chap.\ 3]{KM}. 
Theorem~\ref{prop:FreeV-Hopfian-rk} for the variety of abelian groups can be extended as follows. 

\bc \label{prop:FGAbelian}
\label{abelian}
Let $G$ be a finitely generated abelian group. Then:
$$G\rkS = \left\{
\begin{array}{ll}
G\rkG+1&\mbox{ if $G$ is torsion-free}\\
G\rkG&\mbox{ otherwise}
\end{array}
\right.$$
\ec

\proof
We use the structure theorem for finitely generated abelian groups~\cite[Chap.\ 3]{KM}. 

If $G$ is torsion-free, then $G$ is a free abelian group of finite rank and the claim follows from 
Theorem~\ref{prop:FreeV-Hopfian-rk}\eqref{prop:FreeV-Hopfian-rk-ii}. 

If $G$ is not torsion-free, then $G$ is a finitary direct product of cyclic groups, where at least one of the factor groups is finite. But then we can take a group basis composed by generators of the factor groups, 
which will include an element of finite order, and apply Proposition \ref{torsion}.
\qed

Let $G$ be a group. 
Given $a,b \in G$, as usual write $[a,b] =  aba\inv b\inv$, and 
for $A, B \subseteq G$, 
\[
[A,B] = \langle [a,b] \mid a \in A, \, b \in B \rangle. 
\]
Also, denote the commutator $[G,G]$ of $G$ by $G'$.

\bp \label{prop:FGNilpGp-rkDeriv}
If $G$ is a finitely generated nilpotent group, then 
$G\rkG = (G/G')\rkG$.
\ep

\proof
Assume that $G$ is as in the statement. 
It is clear that 
$(G/G')\rk_G \leq G\rk_G$ and $(G/G')\rk_S \leq G\rk_S$. 
The fact that $G$ is finitely generated and nilpotent implies that 
$G'$ is finitely generated and contained in the Frattini subgroup of~$G$, 
which means that each element of $G'$ can be omitted from every 
set $A$ such that $G = \langle A \rangle$ 
(see~\cite[Th.~2.2.6 and Sec.~16]{KM} and \cite[Chap.~3]{HNeu67}). 
Thus, if ${a_1, \dotsc, a_n \in G}$ are such that $G/G' = \langle a_1G', \dotsc, a_nG'\rangle$, 
then $G = \langle a_1, \dotsc, a_n \rangle$. 
Hence ${G\rk_G \leq (G/G')\rk_G}$, and therefore $G\rk_G = (G/G')\rk_G$. 
\qed

The next result generalizes partially Corollary~\ref{prop:FGAbelian}. 

\bc
Let $G$ be a finitely generated nilpotent group such that $G/G'$ is torsion-free. 
Then 
$G\rkS = G\rkG + 1$. 
\ec

\proof 
By Corollary~\ref{prop:FGAbelian}, the equality 
$(G/G')\rk_S = (G/G')\rk_G + 1$ holds, since $G/G'$ is abelian and torsion-free. 
Then, applying Proposition~\ref{prop:FGNilpGp-rkDeriv}, 
$G\rk_S \leq  G\rk_G + 1 = (G/G')\rk_G + 1 = (G/G')\rk_S \leq G\rk_S$, 
and so $G\rk_S = G\rk_G + 1$.
\qed

An example presented in \cite[31.42]{HNeu67} and attributed to L.~G.~Kov\'acs shows that 
there exists a finitely generated nilpotent group $G$ such that $G$ is torsion-free and $G/G'$ is not. 
In the next example we take that group $G$ and see that $G\rk_S = G\rk_G$. 
Thus the analogue of Corollary~\ref{prop:FGAbelian} cannot be generalized for 
finitely generated nilpotent groups. 
It remains unanswered the question of whether $G\rk_S = G\rk_G$ for any finitely generated nilpotent group $G$ 
such that $G/G'$ is not torsion-free. 

\be \label{ex:Nil-TorsionFree-EqualRanks}
Let $F$ be the relatively free group in the variety of nilpotent groups of class at most~$2$ 
on two generators $a, b$. 
Take the infinite cyclic group $D$ generated by an element $d$ that does not belong to $F$. 
Let $G$ be the free product of $F$ and $D$ amalgamating $[a,b]$ with $d^2$. 
Then $G$ is a finitely generated nilpotent torsion-free group and $G\rkS = G\rkG$.

{\em 
Let us prove these facts. 
Recall, firstly, that a group $H$ being nilpotent of class at most~$2$ means that 
$[H',H] = \{ 1 \}$, which is equivalent to saying that 
every commutator $[x,y]$ of elements of~$H$ commutes with every element of $H$. 
Secondly, that $F \cong  F_2/[F_2',F_2]$ (see \cite[Sec.~16.1]{KM} and \cite[Chap.~3]{HNeu67}). 
Thus, it follows from basic properties on commutators of elements of a group that 
\[
F = \text{Gp}\langle a,b,c  \mid c[b,a], \, [a,c], \, [b,c] \rangle . 
\] 
Then 
\[
G = \text{Gp}\langle a,b,c,d \mid c [b,a], \, [a,c], \, [b,c], \, c\inv d^2 \rangle ,  
\]
and hence 
\[
G/G' = \text{Gp}\langle a,b,d \mid [a,b], \, [a,d], \, [b,d], \, d^2 \rangle . 
\]
It follows that $G/G' \cong \ZZ \times \ZZ \times \ZZ_2$, whence 
$(G/G')\rk_S = (G/G')\rk_G = 3$ by Corollary~\ref{prop:FGAbelian}. 
The group $G$ is nilpotent, since $F$ and $D$ are nilpotent, and is obviously finitely generated.
Then $G\rk_G = 3$ by Proposition~\ref{prop:FGNilpGp-rkDeriv}. 
Since $G = \langle a,b,d \rangle$, by~Proposition~\ref{prop:GgGen-SmgGen} we have 
$G = \bigl\{ a,\, b,\, d,\, a\inv b\inv d\inv\bigr\}^+$. 
However, in $G$, from $d^2 = [a,b]$ we obtain $d = [a,b]d\inv = aba\inv b\inv d\inv$, whence 
$G = \bigl\{ a,\, b,\, a\inv b\inv d\inv\bigr\}^+$. 
Thus $G\rk_S \leq 3$. 
Then 
$3 = (G/G')\rk_G \leq G\rk_G \leq G\rk_S \leq 3$, 
and therefore $G\rk_G = 3 = G\rk_S$.  
Clearly, $G$ is torsion-free and $G/G'$ has elements of order two. 
}
\ee 

The next example shows that $G\rkS = G\rkG+1$ may hold for groups $G$ which are not torsion-free.


\be
\label{ntf}
Let $G$ be the group defined by {\rm Gp}$\la a,b \mid [a,b]^2\ra$. 
Then $G$ is not torsion-free and $G\rkS = G\rkG+1$. 

{\em 
Indeed, by a well-known result of Karrass, Magnus and Solitar (see
\cite[Proposition II.5.17]{LS}), the one-relator group defined by a presentation 
Gp$\la A \mid r^m\ra$, where $A = \{a,b\}$, $r$~is not a proper power in~$\tF_A$ and $m \geq 2$, 
has elements of order $m$, namely $r$. 
Hence $G$ is not torsion-free. 

On the other hand, by sending $a,b$ to a basis of $\tFA_2$, we build a surjective morphism 
${\p \colon G \to \tFA_2}$. 
Since $\tFA_2\rkG = 2$, it follows that $G\rkG = 2$ and so $G\rkS = 3$ by Proposition \ref{hi}.
}
\ee

It is natural to ask if when a group $G$ is finitely generated and 
$G\rkS = G\rkG$, any group basis of $G$ is a semigroup basis. 
The answer is no, independently of $G$ being torsion-free or not, 
as the following example shows.

\be 
(i) Let $G = \ZZ \times \ZZ_2$, which is not a torsion-free group. 
The set $\{ (1,0), (0,1)\}$ is a group basis of $G$, 
and hence $G\rkS = G\rkG$ by Proposition~\ref{torsion}. 
However $\{ (1,0), (0,1)\}$ is not a semigroup basis of~$G$. 


(ii) Let us take the group $G$ of Example~\ref{ex:Nil-TorsionFree-EqualRanks}, 
which is torsion-free. 
We saw that $G\rkS = G\rkG$ and that 
the set $\{ a, b, d \}$ is a group basis of~$G$. 
However $\{ a, b, d \}$ is not a semigroup basis of~$G$. 
\ee

\bp
\label{qp}
Let $G$ be the group defined by the presentation {\rm Gp}$\la A \mid S\ra$, 
where each letter of~$A$ occurs in some word of $S \cap A^+$. Then $G\rkS \leq |A|$.
\ep

\proof
First, notice that in the statement $A^+$ denotes the subsemigroup of~$\tF_A$ generated by~$A$, 
i.e.\ the free semigroup on~$A$. 
Let $\p \colon \tF_A \to G$ be the canonical projection. It suffices to show that $G = (A\p)^+$. 

Let $a \in A$. There exist $u,v \in A^*$ such that $uav \in S$, hence $(uav)\p = 1$. 
Thus 
$$a\inv\p = (a\inv(u\inv(uav)u))\p = (vu)\p = (vu^2av)\p \in (A\p)^+.$$
Since $\tF_A = (A \cup A\inv)^+$, it follows that $G = \tF_A\p = (A\p)^+$ as required.
\qed

Next we look at some groups that arise naturally in Algebraic Topology. 
We start by giving an example of such a group $G$ that is not nilpotent, 
is torsion-free, and where the equality 
$G\rkS = G\rkG$ holds.

\be
\label{kb}
Let $\pi_1(K)$ denote the fundamental group of the Klein bottle $K$. Then $\pi_1(K)$ is torsion-free and $(\pi_1(K))\rkS = (\pi_1(K))\rkG = 2$. 

{\em 
To see this, write $G = \pi_1(K)$. 
It is well-known that $G$ can be defined by the presentation Gp$\la a,b \mid a^2b^2\ra$ 
and is torsion-free noncyclic~\cite{AFR}. 
Hence 
$G\rkS \leq 2$ by Proposition \ref{qp}. On the other hand, $G\rkG = 2$ since $G$ is noncyclic. 
Thus  
$G\rkS = G\rkG = 2$ as claimed.
}
\ee


We conclude with the analysis of our problem in the case of arbitrary surface groups. 
A {\em surface group} is the fundamental group of a connected closed 
(i.e.~compact without boundary) surface  (see~\cite{AFR,Hat02} for details). 

\bt
\label{surf}
Let $M$ be a connected closed surface, and let $\pi_1(M)$ be its fundamental group. 
Then
$$(\pi_1(M))\rkS = \left\{
\begin{array}{ll}
(\pi_1(M))\rkG+1&\mbox{ if $M$ is orientable}\\
(\pi_1(M))\rkG&\mbox{ if $M$ is non-orientable}
\end{array}
\right.$$
\et

\proof
Write $G = \pi_1(M)$. 

Assume first that $M$ is orientable (of genus $g$). From~\cite{AFR}  
the group $G$ is defined by the presentation 
$${\rm Gp}\la a_1,b_1, \ldots, a_g,b_g \mid [a_1,b_1]\dotsm [a_g,b_g]\ra.$$
By sending the generators $a_i,b_i$ to a basis of $\tFA_{2g}$, we build a surjective morphism 
$\p \colon G \to \tFA_{2g}$. 
Since $\tFA_{2g}\rkG = 2g$, it follows that $G\rkG = 2g$, 
and so $G\rkS = 2g+1$ by Proposition \ref{hi}.

Next assume that $M$ is non-orientable (of genus $g$). 
In this case, 
$G$ is defined by the presentation 
$${\rm Gp}\la a_1, \ldots, a_g \mid a_1^2\dotsm a_g^2 \ra.$$
Let $H$ be the direct product of $g$ copies of the cyclic group $C_2$, 
which is defined by the presentation
$${\rm Gp}\la a_1, \ldots, a_g \mid a_1^2, \ldots, a_g^2, \, [a_i,a_j] \, (1 \leq i < j \leq g)\ra.$$
It is straightforward to check that $H\rkG = g$. 
Considering the canonical surjective morphism $\p \colon G \to H$, it follows that $G\rkG = g$. 
In view of Proposition~\ref{qp}, $G\rkS \leq g$. 
Therefore $g = G\rkG \leq G\rkS \leq g$, and so $G\rkS = G\rkG$.
\qed

\section*{Acknowledgements}

The first and second authors were supported by FCT (Portugal)  
through project \linebreak 
UID/MULTI/04621/2013 of CEMAT-CI\^{E}NCIAS. 

The third author was partially supported by CMUP (UID/MAT/00144/2013), which is funded by FCT (Portugal) with national (MEC) and European structural funds (FEDER), under the partnership agreement PT2020.

\bigskip

{\sc M\'ario J.\ J.\ Branco, Departamento de Matem\'atica and CEMAT-CI\^{E}NCIAS, 
Faculdade de Ci\^encias, Universidade de Lisboa, Campo Grande, 1749-016 Lisboa, Portugal} 

{\em E-mail address:} mjbranco@fc.ul.pt

\bigskip

{\sc Gracinda M.\ S.\ Gomes, Departamento de Matem\'atica and CEMAT-CI\^{E}NCIAS, 
Faculdade de Ci\^encias, Universidade de Lisboa, Campo Grande, 1749-016 Lisboa, Portugal}

{\em E-mail address}: gmcunha@fc.ul.pt

\bigskip

{\sc Pedro V.\ Silva, Centro de
Matem\'{a}tica, Faculdade de Ci\^{e}ncias, Universidade do
Porto, R.~Campo Alegre 687, 4169-007 Porto, Portugal}

{\em E-mail address}: pvsilva@fc.up.pt


\begin{thebibliography}{99}

\bibitem{AFR} 
 P.~Ackermann, B.~Fine, and G.~Rosenberger, 
 On surface groups: Motivating examples in combinatorial group theory, 
 in C.~Campbell, M.~Quick, E.~Robertson, and G.~Smith (Eds.), 
  {\em Groups St.~Andrews 2005 (vol.~1)}, 
 London Math.\ Soc.\ Lecture Note Series, vol.~339, p.~96--129. Cambridge University Press, 2007. 


  
%
%

\bibitem{Benois87} M.\ Benois, 
 Descendants of regular language in a class of rewriting systems: algorithm and
 complexity of an automata construction,  
 in {\em Rewriting techniques and applications}, Lecture Notes in Comput.\ Sci., vol.~256, 
 p.~121--132. Springer-Verlag, 1987. 

%
%
%

\bibitem{DEM16} I.~Dolinka, J.~East, and J.~Mitchell, 
 Idempotent rank in the endomoprhism monoid of a nonuniform partition, 
 Bull.\ Aust.\ Math.\ Soc.\ 93, no.\ 1 (2016), 73--91. 

\bibitem{Eps}
D.~B.~A.~Epstein, J.~W.~Cannon, D.~F.~Holt, S.~V.~F.~Levy, M.~S.~Paterson, and W.~P.~Thurston, 
{\em Word processing in groups}, 
Jones and Bartlett Pub., 1992.


\bibitem{FQ11} V.~H.~Fernandes and T.~M.~Quinteiro, 
 On the monoids of transformations that preserve the order and a uniform partition, 
 Comm.~Algebra {\bf 39}, no.\ 8 (2011), 2798--2815.

%
%

\bibitem{Gray14} R.~Gray, 
 The minimal number of generators of a finite semigroup, 
 {\em Semigroup Forum} {\bf 89}, no.\ 1 (2014), 135--154.


%
 
\bibitem{GH87} G.~M.~S.~Gomes and J.~M.~Howie, 
 On the ranks of certain finite semigroups of transformations, 
 Math.\ Proc.\ Cambridge Philos.\ Soc.\ {\bf 101}, no.\ 3 (1987), 395--403.
 
\bibitem{Hat02} 
 A.~Hatcher, 
 {\em Algebraic Topology}, 
 Cambridge University Press, 2002. 
 

%
%

\bibitem{HEO} 
D.~F.~Holt, B.~Eick and E.~A.~O'Brien, 
{\em Handbook of computational group theory}, 
Chapman \& Hall CRC, 2005.

\bibitem{HU79} J.~E.~Hopcroft and J.~D.~Ullman, 
{\em Introduction to Automata Theory, Languages, and Computation}, Addison-Wesley, 1979. 


%



%

\bibitem{KM} M.~I.~Kargapolov and Ju.~I.~Merzjlakov, 
{\em Fundamentals of the Theory of Groups}, Springer-Verlag, 1979. 

 
\bibitem{LS} R.~C.~Lyndon and P.~E.~Schupp, {\em Combinatorial Group
 Theory}, Springer-Verlag, 1977. 
 
%

\bibitem{HNeu67} H.~Neumann, 
 {\em Varieties of Groups}, Springer-Verlag, 1967. 


%


%
%
%


\end{thebibliography}
\end{document}